# Please, Do Tell
Yvonne Lai

"Be a guide on the side, not a sage on the stage." "Ask, don't tell." "Math is not a spectator sport." "Lecturing is educational malpractice." These slogans rally some mathematicians to teach classes that feature "active learning", where lecturing is eschewed for student participation. Yet as much as I believe that students must do math to learn math, I also find blanket statements to be more about bandwagons than considered reflection on teaching. In this column, I urge us to think through the math we offer students and how we set up students to learn. Although I draw primarily from my experiences teaching proofs in abstract algebra and real analysis, the scenarios extend to other topics in first year undergraduate education and beyond.

**Showing and telling the public and private spaces**
If you have ever taught real analysis, you have seen students struggle. The nested quantifiers, new proof structures, and abstraction can feel like an avalanche. Those students who succeed find new mathematical doors to measure theory, topology, and more. Many students don't though, and instead experience real analysis as an exercise in frustration and a message that they cannot do mathematics. The rewards of real analysis may be rich, but most emerge poor.

When I first taught real analysis proofs, such as proving the divergence or convergence of a sequence, I emphasized conceptual understanding and lots of hands-on activities on the definition of a limit. Many students claimed that after these activities, limits made more sense than they did in calculus. Yet they were not able to write proofs. They were unable to connect scratch work with $\epsilon$ and $N$ to writing a clean argument that a sequence converged or diverged. Consequently, some saw no use in the scratch work at all. I encouraged, begged, and pleaded with my students to go through the scratch work, to no avail.

Then I encountered Manya Raman Sundstrom's dissertation, on proof and justification in university calculus, where she used "public" and "private" to distinguish scratch work from the proof one might see in a textbook or perfect problem set write-up.[1] I was inclined initially to dismiss the terminology—after all, "scratch work" suited me just fine, and I heard plenty of others use it—but there was something about the new terms that stuck with me. In any case, my students never seemed to benefit from talk of "scratch work"; something else was needed.

I developed a little talk about "private space" and "public space", likening "private space" to the things you do but don't show, such as trying on different outfits or practicing your smile before a promising job interview or exciting first date. You might talk to friends about the outfits, but you're not going to tell the interviewer or date about it. You might make sure in the privacy of your bathroom that your smile is enthusiastic, but not too eager or desperate; charming but not sleazy; but you're not going to show all these smiles to the interviewer or date. You're just going to show that charming, enthusiastic smile. That's the "public space".

---
[1] Raman, M. J. (2002). Proof and Justification in Collegiate Calculus. Dissertation. University of California-Berkeley. http://www.lettredelapreuve.org/OldPreuve/Newsletter/03Automne/ManyaThesis.pdf

This spiel got students' attention. Even more promising, I noticed that they referenced "private" and "public" voluntarily—an improvement over their de facto embargo of the term "scratch work". But this talk alone didn't help connect the private mathematics, or scratch work, to the final write up, or public mathematics.

I began wondering what it would take to teach the interplay of public and private mathematics. What if I simulated my own private work and its connection to the public proof? I did so. The result was the most success that I have ever had teaching proofs in real analysis. For the first time, I saw students go from doing proofs of $\lim_{n\to\infty} \frac{1}{n+1} = 0, n \in \mathbb{N}$ in class to proving in homework that $\lim_{n\to\infty} \frac{-3n^3 - \frac{2}{5}n + 2}{5n^3 - n} = -\frac{3}{5}, n \in \mathbb{N}$. For those that did not prove the latter, feedback seemed to make sense; when I asked students to revise the proof, most students only needed one or two more tries. On the midterm, the vast majority of students aced a similar limit proof.

Later, on the first day of $\epsilon$-$\delta$ proofs, I once again modeled private work and how it translated to public proof, for a proof of the continuity of $f + g$ given the continuity of $f$ and $g$, and then $-3g$ given the continuity of $g$. I asked students to prove $Af + Bg$ given the continuity of $f$ and $g$ at their seats. The entire class stared at me with looks of, "Well, duh..." When I walked around, I saw many correct proofs written so well they could have appeared in a textbook. Even more remarkably, many proved the continuity of $f^2$ given the continuity of $f$. (If you haven't tried this before, it is a little hairy compared to other standard $\epsilon$-$\delta$ proofs. The week the students handed in this homework, one student protested that they felt cheated because they had condensed upwards of 2 pages of private work to less than ½ a page of public work—and I saw many heads nodding vehemently in concurrence. But, they had done it!) I should mention that at my institution, this is a first course in real analysis, taken by math majors and some math minors—the majority of whom do *not* go on to graduate school. Instead, they find jobs related to insurance and banking, data analysis, software development, or secondary education,.

I made a special effort to show the connection between private and public spaces, by designating one board "private space" and another board "public space". I walked through scratch work, gave a rationale for each step. Then, and crucially, I walked over to the public space, and, as I wrote each step in the final proof, explained where in the private space the step came from. Figure 1 shows the dual board work.

For instance, when proving that $\lim_{n\to\infty} s_n = 0$, where $s_n = \frac{1}{n+1}, n \geq 0$, the claim comes from the problem statement. The statement $N(\epsilon) = \frac{1}{\epsilon} - 1$ is the result of the private space. Showing deductions from $n > N$ will sometimes but not always resemble private computations. Altogether, these inferences show that for all $\epsilon > 0$ there is an $N$ such that $n > N$ implies that $s_n$ is closer than $\epsilon$ to 0, and hence the limit is 0. (One point of discussion I have with students, after they have accepted the proof as valid, is what it means when $\frac{1}{\epsilon} < 1$ and $N(\epsilon) < 0$. In this case, $n > 0$ is always greater than $N(\epsilon)$. This is saying that when $\epsilon$ is large enough, all terms of $s_n$ are closer than $\epsilon$ to the limit 0.)

**Figure 1. Boardwork that shows public and private spaces for a convergence proof.**

After the success of this approach with convergence proofs, I asked students to find divergence proofs without modeling a similar process in class—only to have many turn in muddled work. Once I modeled divergence proofs in this way, the majority were able to write divergence proofs. However, I did not want to have to model every possible kind of proof—an impossible task in any case. And so, I decided to find a time to model the process of proving that something satisfies the negation of a definition. I found this opportunity in the unit on continuity and discontinuity. We rehashed divergence proofs, and discussed how divergence proofs were an example of proving that an object satisfies the negation a definition. Then I asked students to brainstorm how to prove that a function was discontinuous. After we generated a method in class, we talked about how the process of finding a structure for proving discontinuity paralleled that of proving divergence. I assigned both proofs of continuity and discontinuity for homework, and did not model a private-to-public space process for discontinuity proofs. Students generally succeeded at writing both continuity and discontinuity proofs.

My conclusion from this episode is that teaching is not just about finding ways to get students to do more, or building conceptual understanding, or even giving students feedback after they've attempted a proof. These things are important, but they are not enough. Convergence proofs and limit proofs are a new mathematical language, compared to the proofs students are likely to have seen previously. When teaching a new mathematical language, we must be utterly transparent about the process of mathematics from beginning to end. This can include taking class time to walk through how private work leads to public proof, as unnatural as this may seem to do in front of a room of students. The time that one takes to do so can save more time later, as well as open doors to more students about how mathematics works.

**Transparency with proof structures**
One critique of the above show-and-tell of public and private space is that they may reduce proofs to procedures, and therefore further the idea that math is about formulas rather than reasoning. It is a critique I worried about when planning these units and distributing handouts to students on various proofs. I believe that one potential way to counter this effect, and to promote proofs as a genre of communication that makes sense, is discussing where proof structures come from and why they work.

For instance, in real analysis, the structure of a convergence proof comes from the definition of convergence. Because the definition specifies, "For every $\epsilon > 0$, there exists an $N$ such that $n > N$ implies … ", we can think of $N$ as a function of $\epsilon$, and we must find $N(\epsilon)$ that $n > N(\epsilon)$ leads to the desired inequality to show convergence. In practice, when teaching this, I write the definition on the board, and then I write the first few lines and the last few lines of a proof, and ask students to think individually, then share with a partner, why this proof structure would actually show that a sequence is convergent. We then discussed the connection as a whole class, underscoring the point that definitions come with criteria, and the proof is about establishing those criteria. After a similar process with $\epsilon$-$\delta$ proofs, a student—who had failed the course two times before enrolling in my course—practically ran up to me after class, her face beaming. She said, "This is the first time that these proofs have made sense." Three weeks later, she received one of the highest marks on the midterm assessing performance on these proofs.

Another instance of talking through proof structures is when proving results with Bezout's Identity in abstract algebra (that given $a, b, x \in \mathbb{Z}$, if $x = an + bm$ where $n, m \in \mathbb{Z}$, then there exists $k \in \mathbb{Z}$ such that $x = k \cdot \gcd(a,b)$; and further, if $x = k \cdot \gcd(a,b)$, then there exists $n, m \in \mathbb{Z}$ such that $x = an + bm$). Many of these proofs require something that may as well be called "algebraic wishing" — you hope to prove an identity by making an equation work. Previously, I would imitate the style that I had seen my own undergraduate professors take. I would show the proof—possibly with some class participation during computations—and, as if a rabbit in a bow tie had hopped out of a top hat, I'd wax eloquent about the elegant surprises that awaited us in mathematics. But as I saw students struggle with these proofs term after term, I realized that joy was conspicuously absent. In its place was only frustration. I introduced the term "algebraic wishing" and pulled no magic tricks. Instead, I talked about how there are some proofs in math that are about making equations work, and the technique for making equations work is to use known identities and definitions. Since then, students have more success with these proofs, and others later that fall into the same category.

Looking back, I found joy in the "magic tricks" that my professors pulled in my undergraduate classes because the idea of algebraic wishing was one that already made sense to me. My

father—a computer programmer who wished he had been a math major—had explained this technique to me patiently, in middle school and high school, each time it came up. In other words, I already had access to a firm conviction that math could and should make sense, because someone had already been transparent with me about how certain proofs worked. Talk of "magic" was fun only because I knew it wasn't really magic.

**Moving beyond slogans and magic**
Slogans and magic tricks can capture an idea once we are already comfortable with that idea. However, to the uninitiated, they are meaningless. Worse, they don't help someone reason through their own actions to improve their craft. We can find slogans and magic tricks in both mathematics and mathematics teaching.

In mathematics, slogans and magic tricks come in the form of phrases such as "just use the definition", "write down what you know", or "show your work". These phrases each point to a helpful idea, but the phrases alone can't teach the ideas. Instead, students need experiences to help them make sense of these phrases and see how different instantiations of these phrases can fit together in a whole. When instructors attend to students' work and hear how students talk, then, using these observations over time, they can refine the feedback they offer, problems they assign, and explanations they give.

In teaching, slogans and magic tricks come in the form of phrases like "student-centered" and "instructor-centered". These phrases can mean something, but the phrases themselves don't communicate much. Just because students have done math in class or have talked to each other about math doesn't mean that they have learned math. On the other hand, without having students do math in front of us, be it by themselves or with a group, it's hard to tell what they are learning, and whether we should be changing our plans. Some of the most important tools we have, as instructors, are our eyes and ears. Whatever our goals for students are—be it understanding limit proofs or using Bezout's Identity—our students can teach us what they know. When we look and listen, and identify ideas and practices to make more transparent, we can become better teachers.